\documentclass[12pt]{amsart}

\theoremstyle{plain}
\newtheorem{theorem}                 {Theorem}      [section]

\newtheorem{corollary}    [theorem]  {Corollary}
\newtheorem{lemma}        [theorem]  {Lemma}

\theoremstyle{definition}

\newtheorem{remark}       [theorem]  {Remark}
\newtheorem{definition}   [theorem]  {Definition}

\numberwithin{equation}{section}

\def \r{\mbox{${\mathbb R}$}}
\def \s{\mbox{${\mathbb S}$}}

\def \sp #1#2{\langle #1,#2 \rangle}
\def \hwc{horizontally (weakly) conformal}

\def \ha{harmonic map} 

\def \f{\mbox{$\phi$}}
\def \p{\mbox{$\psi$}}
\def \pf{\mbox{$\psi\circ\phi$}}
\def \tpf{\mbox{${\scriptstyle \psi\circ\phi}$}}
\def \tf{\mbox{${\scriptstyle \phi}$}}
\def \tp{\mbox{${\scriptstyle \psi}$}}

\def \tN{\mbox{${\scriptscriptstyle N}$}}

\def \ts{\mbox{${\tiny \s}$}}

\DeclareMathOperator{\trace}{trace}

\DeclareMathOperator{\Index}{index}
\DeclareMathOperator{\nul}{nullity}

\DeclareMathOperator{\id}{Id}

\DeclareMathOperator{\rank}{rank}

\begin{document}

\begin{abstract}
We prove that harmonic morphisms preserve the Jacobi operator along
harmonic maps.
We apply this result to prove infinitesimal and local rigidity (in the
sense of Toth) of harmonic morphisms to a sphere.
\end{abstract}

\title{Harmonic morphisms and the Jacobi operator}

\author{Stefano Montaldo}
\address{Universit\`a degli Studi di Cagliari\\
Dipartimento di Matematica\\
Via Ospedale 72\\
09124 Cagliari}
\email{montaldo@unica.it}
\author{John C. Wood}
\address{Department of Pure Mathematics\\University of Leeds\\
Leeds LS2 9JT \\GB}
\email{j.c.wood@leeds.ac.uk}
\date{October 1999}

\maketitle

\section{Harmonic morphisms}

 \emph{Harmonics maps} ${\phi}:{(M,g)}\to{(N,h)}$ between two smooth
 Riemannian manifolds are critical points of the energy functional
$E(\phi,\Omega)=\frac{1}{2}\int_{\Omega}^{}\, |d\phi|^2\,dv_{g}$ 
for any compact domain $\Omega\subseteq M$ \cite{EelSam64},
i.e. the first variation of the energy vanishes for any smooth
variation of $\phi$. The Euler-Lagrange equation for the energy is the
vanishing of the tension field $\tau_{\phi}=\trace\nabla d\phi$, where
$\nabla$ denotes the connection on
$T^{\ast}M\otimes\phi^{-1}TN$ induced from the Levi-Civita connections
$\nabla^{M}$ on $M$ and
$\nabla^{N}$ on $N$. If $\{e_i\}_{i=1}^{m}$ is a local  orthonormal
frame  on $M$ we have
$\tau_{\phi} =
	\sum_{i=1}^m \left\{\nabla_{e_i}^{\phi}\bigl(d\phi(e_i)\bigr) -
	d\phi\left(\nabla^M_{e_i} e_i\right)\right\}$
where $\nabla^{\phi}$ denotes the pull-back connection on 
$\phi^{-1}TN$.

 Let ${\phi}:{(M,g)}\to{(N,h)}$ be a smooth map between two
Riemannian manifolds.
The tangent space at a point $x\in M$ can be decomposed as
$
T_{x}M=H_{x}\oplus V_{x}
$
where $V_x = \ker(d\phi_x)$ and
$H_x = {V_x}^{\bot}$.
The spaces $V_x$ and $H_x$ are called the \emph{vertical} and
\emph{horizontal space} at the point $x\in M$, respectively.
Let $C_{\tf}=\left\{x\in M\, |\, \rank(d\phi_x)\;\mbox{is not maximal}
\right\}$; points of $C_{\tf}$ are called 
\emph{critical points} of $\phi$.

\begin{definition}
  A map ${\phi}:{(M,g)}\to{(N,h)}$ is called \emph{\hwc\,} if, for every 
  $x\in M$, \emph{either} ${d\phi_{x}}{|_{H_{x}}}$ is conformal 
  and surjective \emph{or} $d\phi_x = 0$.
\end{definition}

In particular, either $x \in M \setminus C_{\tf}$ and $\phi$ is
submersive at $x$ or $x \in C_{\tf}$ and the differential $d\phi_x$ has
rank $0$.
If ${\phi}:{M}\to{N}$ is \hwc, then there exists a function 
${\lambda}:{M\setminus C_{\tf}}\to{\r^{+}}$ such that
$
h_{\tf(x)}(d\phi_{x}(X),d\phi_{x}(Y))\-=\lambda^2 g_{x}(X,Y)\,,
$
for all $X,Y\in H_x$, and $x\in M\setminus C_{\tf}$\,.
The function $\lambda$ can be extended continuously to the whole of $M$ 
by setting $\lambda|_{C_{\tf}}=0$\,;  the extended 
function is called the \emph{dilation function} of $\phi$. Note that
$\lambda^2$ is smooth.

\emph{Harmonic morphisms} between Riemannian manifolds are defined to be
maps which
pull back harmonic functions to harmonic functions; they are
characterized as 
the harmonic maps which are horizontally (weakly) conformal
\cite{Fug78,Ish79A}.  Note that non-constant harmonic morphisms
can only exist if $\dim M \geq \dim N$.  For a bibliography of papers on
harmonic morphisms and an `atlas' of the known examples, see
\cite{hamobib}; for the general theory, see \cite{Baiwoo00}.

\section{Second variation and the Jacobi operator}\label{svfhma}

 Let ${\phi}:{M}\to{N}$ be a \ha, for simplicity we assume
that $M$ is compact.  Then the first variation 
$D_{V}E(\phi)=0$ for all vector fields along $\phi$, where by a {\it
vector field along $\phi$} we mean a smooth section $V\in
\Gamma(\phi^{-1}TN)$ of $\phi^{-1}TN$.
Given a vector field along $\phi$ we consider a
smooth one-parameter variation $\phi_{t}$
($-\epsilon<t<\epsilon$) of $\phi=\phi_{0}$
such that
$
V=\frac{d}{dt}({\phi_{t}})|_{\scriptstyle t =0} 
$.
Let $\nabla^2 = \nabla\circ\nabla^{\phi}$ so that
$$
\nabla^2_{X,Y}W = \nabla^{\phi}_X(\nabla^{\phi}_Y W) -
\nabla^{\phi}_{\nabla^M_X Y} W
	\qquad \bigl(X,Y \in TM, \ W \in \Gamma(\phi^{-1}TN) \bigr) \,;
$$
note that this is tensorial in $X$ and $Y$ but not, in general, symmetric. 
We have the following
second variation formula for the energy  \cite{Maz73,Smi75};
we use Milnor's convention for the curvature as in \cite{EelLem78}.
\begin{eqnarray*}
H_{\tf}(V,V)&=& \frac{d^2 E(\phi_{t})}{d t^2}{\biggr |}_{t=0}=\int_{M}
\sp{\Delta^{\tf} V
- \trace R^{\tN}(d\phi,V)d\phi}{V}\,v_{g}\\
&=& \int_{M} \sp{J^{\tf}V}{V}\, v_{g}\,.
\end{eqnarray*}
Here $\Delta^{\tf}$ is the Laplacian on sections of $\phi^{-1}TN$ given in
a local orthonormal frame $\{e_i\}$ on $M$
by 
$$\Delta^{\tf}=-\trace \nabla^2=-\sum_{i=1}^m \nabla^2_{e_i,e_i} =
 	-\sum_{i=1}^m \bigl\{\nabla^{\tf}_{e_i}\nabla^{\tf}_{e_i}
	-\nabla^{\phi}_{\nabla_{e_i}e_i}\bigr\} \,.
$$

\noindent The operator
  $J^{\tf} = \Delta^{\tf} - \trace R^{\tN}(d\phi,\cdot)d\phi$ is called the
  \emph{Jacobi} operator of $\phi$; the elements of its kernel are
called 
  \emph{Jacobi fields} and form a subspace $J(\phi)\subset
  \Gamma(\phi^{-1}TN)$.

\begin{definition}
  The \emph{index} of $\phi$ is the dimension of the largest subspace of
  $\Gamma(\phi^{-1}TN)$ on which $H_{\phi}$ is negative definite. The {\it
nullity}
  of $\phi$ is the dimension of the kernel of $J^{\tf}$. Alternatively,
  the index of $\phi$ is the sum of the multiplicities of the 
  negative eigenvalues $\lambda$ of $J^{\phi}$ and the  nullity of
$\phi$ is the 
  multiplicity of the eigenvalue $0$ of $J^{\tf}$.
\end{definition}

Note that the Jacobi operator is a linear elliptic self-adjoint operator
with positive principal part $\Delta^{\tf}$.
It follows from standard elliptic
theory (cf. \cite{Maz73}) that the index and the nullity are both finite.

 Our main result is that harmonic morphisms preserve the Jacobi
operator along harmonic maps as follows:

\begin{theorem}\label{teo-jac}
  Let ${\phi}:{M}\to{N}$ be a harmonic morphism, 
  ${\psi}:{N}\to{P}$ a harmonic map and  $V$ a vector 
  field  along $\psi$. Then
  the Jacobi operator for the vector field $V\circ\phi$ 
  along $\psi\circ\phi$ is given by
  $$
  J^{\tpf}(V\circ\phi)=\lambda^2 J^{\tp}(V)\circ\phi\,,
  $$
  where $\lambda$ is the dilation function of $\phi$.
  In particular, if $V$ is a Jacobi field along $\psi$, then 
  $V\circ\phi$ is a Jacobi field along $\psi\circ\phi$.
\end{theorem}
\proof
We first remark that the composition $\pf$ is harmonic. Set $W=V\circ
\phi\in\Gamma(\phi^{-1}\psi^{-1}TP)$.
 By definition of pull-back connection, for any $X\in \Gamma(TM)$, we
have 
\begin{equation}\label{eq-pb}
\nabla^{\tpf}_{X}W=\left\{\nabla^{\psi}_{d\phi(X)}V\right\}\circ\phi.
\end{equation}
Let $\gamma$ be a smooth curve in $M$ tangent to $X$ and extend $X$ to a
vector field along $\gamma$.
Then, by \eqref{eq-pb}, both sides of the next equation are
well-defined and we have
$$
\nabla^{\tpf}_{X}(\nabla^{\tpf}_{X}W)=\left\{\nabla^{\psi}_{d\phi(X)}
 \bigl(\nabla^{\psi}_{d\phi(X)}V\bigr)\right\}\circ\phi \,.
$$
Again, by \eqref{eq-pb},
$$
\nabla^{\tpf}_{\nabla^M_{X}X}W=\left\{\nabla^{\psi}_{d\phi(\nabla^M_{X}X)}V\right\}\circ\phi
$$
so that
\begin{eqnarray*}
\nabla^2_{X,X}W&=
&\nabla^{\tpf}_{X}\bigl(\nabla^{\tpf}_{X}W\bigr)-\nabla^{\tpf}_{\nabla_{X}X}W\\
&=&\left\{\nabla^{\tp}_{d\phi(X)}\bigl(\nabla^{\tp}_{d\phi(X)}V\bigr)-\nabla^{\psi}_{d\phi(\nabla^M_{X}X)}V\right\}\circ\phi
\,.
\end{eqnarray*}
Since $\phi$ is harmonic, for any orthonormal frame $\{e_i\}$ on $M$ we
have
$$
\sum_{i=1}^m d\phi\bigl(\nabla^M_{e_i} e_i\bigr)=
\sum_{i=1}^m \nabla^N_{d\phi(e_i)} d\phi(e_i)
$$
so that
\begin{eqnarray*}
\trace \nabla^2 W&=&
\sum_{i=1}^m\bigl\{\nabla^{\tp}_{d\phi(e_i)}\nabla^{\tp}_{d\phi(e_i)}V-
\nabla^{\psi}_{\nabla^N_{d\phi(e_i)}d\phi(e_i)}V\bigr\}\circ\phi\\
&=&
\sum_{i=1}^m \bigl\{\nabla^2_{d\phi(e_i),d\phi(e_i)}V\bigr\}\circ\phi \,.
\end{eqnarray*}
Now choosing the bases $\{e_i\}$, $\{f_i\}$ at the points $x$,
$\phi(x)$ such that $d\phi(e_i)=\lambda f_{i}$
for $i=1,\ldots,n$ and $d\phi(e_i)=0$ otherwise, we have
\begin{eqnarray*}
\trace \nabla^2 W&=&\sum_{i=1}^n
\nabla^2_{d\phi(e_i),d\phi(e_i)}V\circ\phi\\
&=&\lambda^2 \trace \nabla^2V \circ\phi.
\end{eqnarray*}
Next, an easy calculation gives
$$
\trace R^P(d(\pf),W)d(\pf)=\lambda^2\trace R^P(d\psi,V)d\psi
$$
so that 
\begin{eqnarray*}
J^{\tpf}(W)&=&J^{\tpf}(V\circ\phi)=\lambda^2\trace\bigl\{-\nabla^2 V -
R^P(d\psi,V)d\psi \bigr\}\circ\phi\\
&=&\lambda^2 J^{\tp}(V)\circ\phi \,.
\end{eqnarray*}
\endproof

\begin{remark}
A version of Theorem~\ref{teo-jac} with the assumption that
the fibres of $\phi$ are minimal was proved by the first author in
\cite{Mon99}.
\end{remark}

\begin{corollary}\label{teo-inst}
  Let ${\phi}:{M}\to{N}$ be a non-constant harmonic morphism
between compact manifolds and let 
  ${\p}:{N}\to{P}$ be a harmonic map. Then 
  \begin{itemize}
  \item[(a)]
  \begin{enumerate} 
  \item[(i)] $\Index(\pf)\geq\Index(\p)$;
  \item[(ii)] $\Index(\phi)\geq\Index(\id^{N})$, in particular, 
  if $\id^{N}$ is unstable, then so is $\phi$;
  \end{enumerate}
\smallskip
  \item[(b)]
  \begin{enumerate} 
  \item[(iii)] $\nul(\pf)\geq\nul(\p)$;
  \item[(iv)] $\nul(\phi)\geq\nul(\id^{N})$.
  \end{enumerate}
  \end{itemize}
\end{corollary}
\proof
(a) For part (i), let $V\in \Gamma(\p^{-1}TP)$ be an
eigenvector with negative eigenvalue, thus
\begin{equation}\label{eig}
J^{\tp}(V)=\alpha V
\end{equation}
for some constant $\alpha <0$. Set $W=V \circ
\phi\in\Gamma\bigl((\pf)^{-1}TP\bigr)$. Then,
by Theorem~\ref{teo-jac},
\begin{eqnarray*}
J^{\tpf}(W)&=& \lambda^2 J^{\tp}(V)\circ\phi \\
&=& \lambda^2 \alpha\, (V\circ \phi) = \lambda^2 \alpha W \,.
\end{eqnarray*}
Now $W$ cannot be identically zero; if it were, $V$ would be zero
on $\phi(M)$; this is an open set since any non-constant
harmonic morphism is an open mapping \cite{Fug78}. 
But then, by the unique continuation theorem \cite{Aro57,Cor56}
applied to \eqref{eig},
$V$ would be identically zero on $N$.  Further, by unique
continuation for harmonic morphisms \cite{Fug78}, $\lambda$ cannot
be zero on an open set. Hence
$$
H_{\phi}(W,W) = \int_M \lambda^2 \alpha \, \sp{W}{W}  < 0 \,.
$$
Now let $V_1,\ldots,V_s$ be linearly independent eigenvectors with negative
eigenvalues.  Then, again by unique continuation,
$W_1=V_1\circ\phi,\ldots,W_s=V_s\circ\phi$ are linearly
independent.
Further, if $V_i$ and $V_j$ correspond to different eigenvalues
$\alpha_i$ and $\alpha_j$,
$$
\alpha_i \int_M \lambda^2 \sp{W_i}{W_j}
	= \int_M \sp{J^{\psi\circ\phi}W_i}{W_j}
	= \int_M \sp{W_i}{J^{\psi\circ\phi}W_j}
	= \alpha_j \int_M \lambda^2 \sp{W_i}{W_j}
$$
so that $H_{\psi\circ\phi}(W_i,W_j) = 0$.
It follows that $H_{\phi}$ is negative definite on the span of the
$W_i$ and the estimate on the index follows.
Part (ii) follows by putting $\psi=\id^{N}$.

(b) Let $V$ be a Jacobi vector field along $\psi$. Then, by 
Theorem~\ref{teo-jac}, so is $W=V\circ\phi$ and we argue as in
part (a). 
\endproof

\section{Toth rigidity of harmonic morphisms}
 Let $\phi:M\to N$ be a harmonic map.
A section $V\in \Gamma(\phi^{-1}TN)$ is said to be a \emph{harmonic
variation\/}
of $\phi$ if  $\phi_{t}=\exp(tV):M\to N$ is harmonic for 
all $t\in\r$.
Let $H(\phi)\subset \Gamma(\phi^{-1}TN)$ denote the set of all harmonic
variations of a given map $\phi:M\to N$.

 Note that, if $\phi_{t}$ is a variation of $\phi$ through harmonic
maps,
then $V=\frac{\partial \tf_t}{\partial t}|_{t=0}$ is a Jacobi field 
along \f\ \cite{EelLem78}. The converse is not always true,
in fact there are examples of Jacobi fields along a
harmonic map which do not arise as variations through harmonic
maps (see, e.g., \cite{Smi72}). Nevertheless we have
the following.

\begin{theorem}[\cite{Tot82B}]\label{teo-toth}
  Let $\f:(M,g)\to N=N(c)$ be a harmonic map from a compact
  Riemannian manifold to a real space
  form of constant curvature $c\neq 0$.
  Then a section $V\in \Gamma(\f^{-1}TN)$ is a 
  harmonic variation if and only if $V$ is a Jacobi field 
  with $\|V\|$ constant and $\trace\sp{d\f}{\nabla^{\phi} V}=0$.
\end{theorem}

 Let  $\f:M\to \s^{n}$ be a harmonic map to the Euclidean
sphere of dimension $n$ ($n \in \{1,2,\ldots\}$) and set
$$
K(\f)=\{V\in J(\f)\,:\,\trace\sp{d\f}{\nabla^{\phi} V}=0\}.
$$

\noindent Then $K(\f)\subset J(\f)$ is a linear subspace and, by Theorem
\ref{teo-toth},
$$
H(\f)=\{V\in K(\f)\,:\, \|V\| \; \mbox{is constant}\} \,.
$$

Recall that a variation $V\in \Gamma(\f^{-1}TN)$
 is called \emph{projectable} if $\f(x)=\f(x')$ implies that
$V(x)=V(x')$.

\begin{definition}[\cite{Tot82B}]
  A harmonic map $\f:M\to\s^n$ is said to be
  \begin{enumerate}
  \item \emph{infinitesimally rigid\/} if, for every projectable
  $V\in K(\f)$, there exists $X\in so(n+1)$ such that the equation 
  $V=X\circ\f$ holds;
  \item \emph{locally rigid\/} if, for every projectable harmonic
   variation $V$, there exists a $1$-parameter subgroup
  $(g_t)\subset O(n+1)$ such that 
  $\f_t=\exp(tV)=g_t\circ\f$ for all $t\in\r$.
  \end{enumerate}
\end{definition}

\begin{lemma}[\cite{Tot82B}]\label{lem-toth}
  Let $\id^{\ts^n}:\s^n\to\s^n$ be the identity map. Then 
  $$K(\id^{\ts^n})=so(n+1).$$
\end{lemma}

 We extend a theorem of Toth \cite{Tot82B} on rigidity of harmonic
Riemannian submersions to submersive harmonic morphisms.

\begin{theorem}

  (i) Any surjective and submersive harmonic morphism $\f:M\to \s^n$ is
infinitesimally rigid.

 (ii) Let $n$ be odd.  Then any surjective and submersive harmonic morphism $\f:M\to \s^n$ is
   locally rigid.
  \end{theorem} 
\proof 
(i) Let $V\in K(\f)$ be projectable. Then since $\phi$ is surjective and
submersive we have $V=X\circ\f$ for some 
$X\in\Gamma(T\s^n)$ and so, from Theorem~\ref{teo-jac}, 
\begin{equation}\label{eq-A}
 J^{\tf}(V)=\lambda^2\,J^{\ts^n}(X)\circ\f\,,
\end{equation}
where $J^{\ts^n}$ denotes the Jacobi operator along 
$\id^{\ts^n}$\,.
Moreover, choosing orthonormal bases $\{e_i\}$ at a point $x\in M$  
and $\{f_i\}$ at $\phi(x)\in N$
such that $d\phi(e_i)=\lambda f_{i}\circ\phi$ for $i=1,\ldots,n$ and
$d\phi(e_i)=0$ otherwise, we have
\begin{eqnarray}\label{eq-B}
 0=\trace\sp{d\f}{\nabla^{\phi} V}&=&\sum_{i=1}^m
\sp{d\phi(e_i)}{\nabla^{\phi}_{e_i}( X\circ\phi)}\\
&=& \lambda^2 \sum_{i=1}^n \sp{f_i}{\nabla_{f_i} X}\circ\phi\nonumber\\
&=&\lambda^2\, \trace\sp{d(\id^{\ts^n})}{\nabla X}\circ\phi\nonumber
\end{eqnarray}
where $\nabla$ denoted the Levi-Civita connection on $S^n$.
Equations \eqref{eq-A} and \eqref{eq-B} imply that $X\in K(\id^{\ts^n})$
and by Lemma \ref{lem-toth}, $X\in so(n+1)$. This proves that \f\ is
infinitesimally rigid.

\noindent (ii) Suppose that $V\in H(\f)$ so that $\|V\|$ is constant. As
in part (i), $V=X\circ\f$ for some $X\in\Gamma(TS^n)$.  Since $\phi$ is
surjective, $\|X\|$ is also constant and, since, as in the first part of
the proof, $X$ is Killing, for any
vector field 
$Y$ on $\s^n$ we have 
$$
\sp{Y}{\nabla_{X}X}=-\sp{X}{\nabla_{Y}X}=-\frac{1}{2}Y\|X\|^2=0.
$$
It follows that $\nabla_{X}X=0$ on $\s^n$ which implies that every
integral curve $t\to g_t(x)$ \ ($t\in\r$, $x\in\s^n$) of $X$ is a
geodesic.
Hence
$$
\f_{t}(x)=\exp(tV_{x})=\exp(tX_{\tf(x)})=g_t(\f(x))
$$
for all $x\in M$, and, since $X\in so(n+1)$, we have $(g_t)\subset O(n+1)$.

\endproof

 We remark that part (ii) is true for $n$ even but says nothing new as there
are no non-trivial projectable harmonic variations $V \in K(\phi)$ in that case.


\end{document}